\documentclass[final,leqno,onefignum,onetabnum]{siamltex1213}

\usepackage{amsfonts}
\usepackage{amsmath}
\usepackage{amssymb}
\usepackage{epsfig}
\usepackage{capt-of}
\usepackage{multirow}

\def\brdf{{\beta_{RDF}}}

\newcommand{\csb}[1]{\color{black}#1\normalcolor}
\title{Accelerating the Uzawa Algorithm} 

\author{Nguyenho~Ho\thanks{Department of Mathematical Sciences,
    Worcester Polytechnic Institute, Worcester, MA 01609-2280, USA
    (\email{nho@wpi.edu, sdolson@wpi.edu, walker@wpi.edu}). The first
    author is the corresponding author. This work was supported in
    part by US National Science Foundation grants DMS~1122461 and
    DMS~1337943, and US Department of Energy award number DE-SC0004880.}  
\and Sarah~D.~Olson\footnotemark[1] 
\and Homer~F.~Walker\footnotemark[1] 
}

\begin{document}
\maketitle

% \footnotetext[2]{Support from DOE \#}
% \footnotetext[3]{Department of Mathematical Sciences, Worcester Polytechnic Institute, Worcester, MA 01609-2280, USA }
% \footnotetext[4]{Support from NSF DMS 1122461.}

\slugger{sisc}{xxxx}{xx}{x}{x--x}%slugger should be set to mms, siap, sicomp, sicon, sidma, sima, simax, sinum, siopt, sisc, or sirev

\begin{abstract}
  The Uzawa algorithm is an iterative method for the solution of
  saddle-point problems, which arise in many applications, including
  fluid dynamics. Viewing the Uzawa algorithm as a fixed-point
  iteration, we explore the use of Anderson acceleration (also known
  as Anderson mixing) to improve the convergence. We compare the
  performance of the preconditioned Uzawa algorithm with and without
  acceleration on several steady Stokes and Oseen problems for
  incompressible flows. For perspective, we include in our comparison
  GMRES with two different preconditioners. The results indicate that
  the accelerated preconditioned Uzawa algorithm converges
  significantly faster than the algorithm without acceleration and is
  competitive with the other methods considered.
\end{abstract}


%\begin{AMS}\end{AMS}

\pagestyle{myheadings}
\thispagestyle{plain}
\markboth{Ho, Olson, and Walker}{Accelerating Uzawa}

\section{Introduction}
A saddle-point problem is a block 2$\times$2 linear system of the form
\begin{eqnarray}
\underbrace{\left(\begin{array}{cc}
{A} & {B}^T  \\
{B}& C \end{array} \right)}_{\mathcal{A}}
\left(\begin{array}{c}
{ u}  \\
p \end{array} \right)=
\underbrace{\left(\begin{array}{c}
{ f} \\
g\end{array} \right)}_b .\label{Saddle}
\end{eqnarray}
Here, $A$ is assumed to be invertible but may be symmetric or
non-symmetric.  Applications in which saddle-point problems arise
include computational fluid dynamics \cite{Glowinski84, Quarteroni94},
constrained optimization \cite{Wright92}, linear elasticity,
economics, and many other areas \cite{Benzi05}.

Many methods for solving saddle-point problems have been developed and
analyzed, including direct and iterative methods (see \cite{Benzi05}
for an extensive review). Here, we consider the {\em Uzawa algorithm}
\cite{Uzawa58}, a well-known and easily implemented iterative
method. Extensions of the classical Uzawa algorithm include the
augmented Lagrangian formulation \cite{Benzi11d,Liu01} and the inexact
Uzawa algorithm (also known as the Arrow--Hurwicz method)
\cite{Bramble97,Cao03,Elman94,Elman96,Hu06,Lin06,Queck89,Zulehner01}.
A relaxation parameter appears in the algorithm, and methods for
determining optimal values of this are outlined in \cite{Bacuta11}.
Additionally, the use of preconditioners to increase the speed of
convergence of these methods has been studied by a number of authors 
\cite{Bramble97,Bramble99,Cao03,Cui04,Elman94,Elman96,Klawonn98}.

In practice, the iterates produced by the Uzawa algorithm often
converge undesirably slowly. In this paper, we view the algorithm as a
fixed-point iteration and consider the use of {\em Anderson acceleration}
\cite{Anderson65}, also known as {\em Anderson mixing}, to improve the
convergence. This acceleration method (denoted henceforth by AA) has
enjoyed considerable success in accelerating fixed-point iterations in
electronic-structure computations (see \cite{Fang09} and the
references therein), simulations of transport phenomena
\cite{Calef13,Kuzmin12}, and coupled fluid-solid interface problems
\cite{Ganine13}. It has also been effective in accelerating Picard
iterations in variably saturated flow simulations \cite{Lott12} and
advection-diffusion problems \cite{Lipnikov13}. To the best of our
knowledge, it has not been previously used with the Uzawa algorithm or
other iterations for saddle-point problems.

It is unusual to apply AA in the context of a linear problem such as
\eqref{Saddle}, for which there are many possible alternative solution
methods. Indeed, as we note below, AA coupled with the Uzawa algorithm
is very closely related to GMRES \cite{SaadSchultz86} with a certain
preconditioning applied to \eqref{Saddle}. Our goal is to present AA
as an easily implemented and economical way to mitigate the often slow
convergence of Uzawa iterates and to demonstrate the performance of
the accelerated algorithm on several model problems.
  
In the following, we focus on applying AA to the Uzawa algorithm for
the Stokes and Oseen problems of steady incompressible flow. The use
of AA with any iterative method in this fluids setting appears to be
new. In this context, $C=0$ in \eqref{Saddle} when applying a stable
finite-element discretization. (The extension to non-zero $C$ is
straightforward.)  We briefly consider the classical standard form of
the Uzawa algorithm but focus primarily on the preconditioned form,
since the standard form without preconditioning converges too slowly
to be practical on the problems of interest. In Section~2, we outline
the standard and preconditioned algorithms. In Section~3, we first
note a particular way in which the Uzawa algorithm can be recast as a
fixed-point iteration; we then describe AA and discuss its
relationship with GMRES in this case. In Section~4, we report on
numerical experiments, in which we compare the performance of the
Uzawa algorithm, with and without acceleration, in several Stokes and
Oseen flow scenarios. For perspective, we include the closely related
preconditioned GMRES method and also, in most tests, GMRES
preconditioned with the Relaxed Dimensional Factorization (RDF)
preconditioner developed in \cite{Benzi11b}. In Section~5, we offer a
summary discussion and conclusions.

\section{Variations of the Uzawa algorithm}
\label{varUzawa}
\subsection{Standard Uzawa}\label{standardU}
The standard Uzawa algorithm consists of a coupled iteration for
the variables $p$ and $u$. Given an initial guess for $p_0$, the
iteration is
\begin{subequations}
\begin{align}
u_{k+1}&=A^{-1}(f-B^Tp_k), \label{U1}\\
p_{k+1}&=p_k+\omega( Bu_{k+1}-g), \label{U2}
\end{align}
\end{subequations}
where $\omega>0$ is a relaxation parameter. 
It has been previously shown that if one eliminates $u_{k+1}$ from
\eqref{U2} by using \eqref{U1}, then the Uzawa method is equivalent to
a Richardson iteration, which can be used to determine an optimal
$\omega$ in terms of the maximum and minimum eigenvalues of the Schur
complement $\mathcal{S}=BA^{-1}B^T$ \cite{Elman94}. Bounds obtained in
\cite{Elman94} suggest that the iterates will converge slowly if the
Schur complement is not well conditioned. Several variations of this
standard algorithm have been developed, including preconditioned
versions that we will describe next.

\subsection{{Preconditioned} Uzawa}\label{precondU}
We consider the preconditioned Uzawa algorithm in the general form
given in \cite{Bramble97}:
\begin{subequations}
\begin{align}
{u_{k+1}}&=u_k+Q_A^{-1}\left({f-(Au_k+B^Tp_k})\right),\label{PU1}\\
p_{k+1}&=p_{k}+\omega Q_B^{-1}(Bu_{k+1}-g),\label{PU2}
\end{align}
\end{subequations}
where $Q_A$ and $Q_B$ are preconditioners.\footnote{This iteration is
  called an {\em inexact} Uzawa algorithm in
  \cite{Bramble97,Bramble99}, since the ultimate interest there is in
  iterative solution of the linear subproblems. It is not the same as
  the inexact Uzawa algorithm of \cite{Elman94}.}
When $Q_A=A$, this becomes the preconditioned Uzawa algorithm
considered in \cite{Elman94}; when $Q_B=I$ as well, it reduces to the
standard Uzawa algorithm in \eqref{U1}$-$\eqref{U2}. Note also that
\eqref{PU1}$-$\eqref{PU2} require an initial guess for $u_0$ as well as
$p_0$ if $Q_A\ne A$.

Often, $A$ is symmetric positive-definite (SPD). In this case,
suitable preconditioners $Q_A$ and $Q_B$ should be SPD as well, and
various choices have appeared in the literature. For example,
\cite{Elman94} includes an instance of \eqref{PU1}$-$\eqref{PU2} with
$Q_A=A$ and $Q_B$ a scaled diagonal or tridiagonal matrix derived from
the mass matrix of a finite-element discretization. Other examples
appear in \cite{Bramble97}, with $Q_A$ a multigrid V-cycle and
$Q_B=I$, and in \cite{Cao03}, with $Q_A$ an incomplete-Cholesky
factorization of $A$ and $Q_B$ a scaled identity matrix.  Choices of
$Q_A$ and $Q_B$ when $A$ is non-symmetric are similarly varied. For
example, in both \cite{Bramble99} and \cite{Cui04}, the symmetric part
$A_S={{1}\over{2}}(A+A^T)$ is assumed to be SPD. Choices of $Q_A$ and
$Q_B$ in \cite{Bramble99} include taking $Q_A$ to be a multigrid
V-cycle on $A_S$ and $Q_B$ to be a scaled identity matrix; choices in
\cite{Cui04} include $Q_A=A$ and $Q_B=B(A_S)^{-1}B^T$.

\section{Accelerating preconditioned Uzawa}
\label{accUzawa}
Although the developments in this section are for the preconditioned
Uzawa algorithm \eqref{PU1}--\eqref{PU2}, it is straightforward to
adapt them to the standard algorithm \eqref{U1}--\eqref{U2} by taking
$Q_A=A$ and $Q_B=I$.

\subsection{{Preconditioned} Uzawa as a fixed-point
  iteration}\label{inexact}

We begin by rewriting the preconditioned Uzawa algorithm as a
fixed-point iteration. From \eqref{PU1}--\eqref{PU2}, we have
\begin{equation}
\underbrace{\left( \begin{array}{cc}Q_A & 0 \\ B &
      -\frac{1}{\omega}Q_B \end{array}\right)}_{\mathcal{M}}\left( \begin{array}{c}u_{k+1}\\
    p_{k+1} \end{array}\right)=\underbrace{\left( \begin{array}{cc}
      Q_A-A & -B^T \\ 0 &
      -\frac{1}{\omega}Q_B \end{array}\right)}_{\mathcal{N}}\left( \begin{array}{c}u_k
    \\ p_k \end{array}\right)+\left( \begin{array}{c}f \\
    g \end{array}\right),\label{UzawaMN} 
\end{equation}
which immediately gives the following fixed-point version of preconditioned Uzawa: 
\begin{equation}
  \left( \begin{array}{c}u_{k+1}\\
      p_{k+1} \end{array}\right)=G\left( \begin{array}{c}u_k \\
      p_k \end{array}\right)\equiv
  \mathcal{M}^{-1}\left(\mathcal{N}\left( \begin{array}{c}u_k \\
        p_k \end{array}\right)+\left( \begin{array}{c}f \\ 
        g \end{array}\right)\right). 
\label{UF}
\end{equation}
The explicit form of the right hand side of \eqref{UF} is
\begin{equation*}
\left( \begin{array}{cc}Q^{-1}_A(Q_A-A) & -Q_A^{-1}B^T \\[4pt] \omega
    Q_B^{-1}B Q_A^{-1}(Q_A-A)& I- \omega
    Q_B^{-1}BQ^{-1}_AB^T \end{array}\right)\left( \begin{array}{c}u_k
    \\ p_k \end{array}\right)+\left( \begin{array}{c}Q_A^{-1}f \\[4pt]
    \omega Q_B^{-1}(BQ_{A}^{-1}f-g) \end{array}\right),
\end{equation*}
where $BQ_A^{-1}B^T$ is an approximation of the Schur complement $BA^{-1}B^T$. 
In the case $Q_A=A$, the form of \eqref{UF} is
\begin{equation}
  \left( \begin{array}{c}u_{k+1}\\
      p_{k+1} \end{array}\right)=\left( \begin{array}{cc}0 &
      -A^{-1}B^T \\[2pt] 0 & I-\omega
      Q_B^{-1}BA^{-1}B^T \end{array}\right)\left( \begin{array}{c}u_k
      \\ p_k \end{array}\right)+\left( \begin{array}{c}A^{-1}f \\[2pt]
      \omega
      Q_B^{-1}(B{A}^{-1}f-g) \end{array}\right).\label{StokesU}
\end{equation}
Formally, the iteration \eqref{StokesU}, like \eqref{UF}, requires an
initial guess for $u_0$ as well as $p_0$. However, in \eqref{StokesU},
this can be arbitrary. 

\subsection{Anderson acceleration}\label{Anderson}

The outline of AA given below applies to a general fixed-point
iteration $\xi_{i+1}=G(\xi_k)$ that begins with some initial value
$\xi_o$. In the context of interest, $\xi=(u,\hspace{.1cm} p)^T$ and
$G$ is defined in \eqref{UF}. In the algorithm, $m$ denotes the
maximum number of the stored residuals, which is necessarily finite in
practice. 
% with initial value $\xi_o$ and $m\geq 1$. We set $\xi_1=G(\xi_o)$
% where $G$ is defined in \eqref{UF}. The Anderson Acceleration
% algorithm \cite{Walker11} is as follows:
\begin{eqnarray*}
&&\text{\bf Anderson Acceleration}\\
&& \text{Given $\xi_0$ and $m\ge1$, set $\xi_1=G(\xi_0)$.}\\
&&\text{For } k=1,2,...\\
&&\hspace{10mm}\text{ Set } m_k=\min\{m,k\}.\\
% &&\hspace{10mm}\text{ Set } \tilde{F}_k=(\tilde{f}_{k-m_k},\cdots,
% \tilde{f}_k), \text{ where } \tilde{f}_i=G(\xi_i)-\xi_i\\ 
&&\hspace{10mm}\text{ Set } F_k=(\tilde{f}_{k-m_k},\cdots,
\tilde{f}_k), \text{ where } \tilde{f}_i=G(\xi_i)-\xi_i.\\ 
&&\hspace{10mm}\text{ Determine }
\alpha^{(k)}=\left(\alpha^{(k)}_0,\cdots, \alpha^{(k)}_{m_k}\right)^T
\text{ that solves }\\ 
% &&\hspace{20mm}\min_{\alpha=(\alpha_0,\cdots,
%   \alpha_{m_k})^T}\|\tilde{F}_k\alpha\|_2 \hspace{5mm}\text{
%   s.t. }\hspace{5mm} \sum_{i=0}^{m_k}\alpha=1\\ 
&&\hspace{20mm}\min_{\alpha=(\alpha_0,\cdots,
  \alpha_{m_k})^T}\|F_k\alpha\|_2 \hspace{5mm}\text{
  s.t. }\hspace{5mm} \sum_{i=0}^{m_k}\alpha=1.\\ 
&&\hspace{10mm}\text{ Set }
\xi_{k+1}=\sum_{i=0}^{m_k}\alpha_i^{(k)}G(\xi_{k-m_k+i}).
\end{eqnarray*}

The rationale underlying AA is that, on a linear fixed-point problem
such as the problem of interest here, one has
$\sum_{i=0}^{m_k}\alpha_i^{(k)}G(\xi_{k-m_k+i})=G(\sum_{i=0}^{m_k}\alpha_i^{(k)}\xi_{k-m_k+i})$. 
It follows that, in this case, $\xi_{k+1}$ is obtained by applying $G$
to the point within the affine subspace containing $\xi_k$, \dots,
$\xi_{k-m_k}$ that has minimal fixed-point residual.

In practice, the constrained least-squares problem in $m_k+1$
variables is usually reformulated as an equivalent unconstrained
problem in $m_k$ variables. See \cite{Fang09,Walker11} for a
particular reformulation that can be implemented efficiently and, in
our experience, has been numerically sound. In any case, as the
algorithm proceeds, both the number of stored residual vectors and the
amount of arithmetic per iteration increase up to maxima determined by
$m$. It is important to note, though, that only a single
$G$-evaluation is required at each iteration, as in the underlying
fixed-point iteration. Thus the cost of applying AA to a fixed-point
iteration is only some additional storage and arithmetic (usually
modest in practice); no additional $G$-evaluations are required.

\subsection{Matrix-splitting preconditioners}\label{PGMRES}

As noted in the introduction, the Uzawa algorithm accelerated with AA
is closely related to GMRES applied to \eqref{Saddle} with a certain
preconditioning. To describe this relationship in more detail, we note
that the matrix $\mathcal{A}$ in \eqref{Saddle} with $C=0$ can be split as
follows:
\begin{equation}
\mathcal{A}=\mathcal{M}-\mathcal{N}=\left(\begin{array}{cc}Q_A & 0 \\ B & -\frac{1}{\omega}Q_B \end{array}\right)-\left(\begin{array}{cc}Q_A-A & -B^T \\ 0 & -\frac{1}{\omega}Q_B \end{array}\right),\label{Split}
\end{equation}
where $\mathcal{M}$ and $\mathcal{N}$ are defined in \eqref{UzawaMN}
for $Q_A$ and $Q_B$ in \eqref{PU1}--\eqref{PU2}.\footnote{When
  $Q_A=A$, this splitting has been considered by others; see
  Section~8.1 of \cite{Benzi05} and the references therein.}

The fixed-point form of preconditioned Uzawa \eqref{UF} is a
stationary iteration determined by this splitting. Suppose that AA is
applied to this stationary iteration without truncation, i.e., with
$m_k=k$ at each iteration, so that all previous residuals are used in
the least-squares problems. Suppose also that unrestarted GMRES is
applied to the left-preconditioned system
\begin{equation}
\mathcal{M}^{-1}\mathcal{A}\left(
\begin{array}{c} u\\p\end{array}\right)=
\mathcal{M}^{-1} \left( \begin{array}{c}f\\ g \end{array}\right), 
\label{LPSystem}
\end{equation}
starting from the same initial approximate solution. Denote the $k^{th}$
iterates of the two algorithms by $\xi_k^{\rm AA}$ and $\xi_k^{\rm
  GMRES} $, respectively, and denote the $k^{th}$ GMRES residual by
$r_k^{\rm GMRES}$. Then, provided GMRES does not stagnate before the
solution is found, Corollary~2.10 in \cite{Walker11} asserts that the
two algorithms are ``essentially equivalent'' in the sense that, until
the solution is found, $\xi_k^{\rm GMRES} =
\sum_{i=0}^{m_k}\alpha_i^{(k)}\xi_i^{\rm AA}$ and $\xi_{k+1}^{\rm AA}
= G(\xi_k^{\rm GMRES})$, where $G$ is the fixed-point map in
\eqref{UF}.

This ``essential equivalence'' seldom strictly holds in practice,
since it is almost always necessary to truncate AA and to restart
GMRES. We explore in Section~4 how well this ``essential equivalence''
is borne out in the problems of interest there. For the present, we
note that 
\begin{equation*}
  \mathcal{M}^{-1} =\left( \begin{array}{cc}Q_A^{-1} & 0 \\[4pt]
\omega Q_B^{-1}BQ_A^{-1} & -\omega Q_B^{-1} \end{array}\right), 
\end{equation*}
and so the main requirement of using this preconditioner in GMRES is
one solve with each of $Q_A$ and $Q_B$ per iteration. Additionally,
the preconditioned system \eqref{LPSystem} corresponds to
\begin{equation*}
  \left(\begin{array}{cc}Q_A^{-1}A & Q_A^{-1}B^T \\[4pt] \omega
      Q_B^{-1}BQ_{A}^{-1}(A-Q_{A}) & \omega
      Q_{B}^{-1}BQ_{A}^{-1}B^T \end{array}\right)\left( \begin{array}{c}u\\
      p \end{array}\right)=\left( \begin{array}{c} Q_A^{-1}f \\[4pt] \omega
      Q_B^{-1}(B{Q_A}^{-1}f-g) \end{array}\right),\label{SplitDetail} 
\end{equation*}
When $Q_A=A$, this simplifies to
%\eqref{SplitDetail} simplifies to
\begin{equation}
  \left(\begin{array}{cc}I & A^{-1}B^T \\[2pt] 0 & \omega
      Q_{B}^{-1}B{A}^{-1}B^T \end{array}\right)\left( \begin{array}{c}u\\
      p \end{array}\right)=\left( \begin{array}{c}A^{-1}f \\[2pt] \omega
      Q_B^{-1}(B{A}^{-1}f-g) \end{array}\right).\label{PBigSchur} 
\end{equation} 

We refer to GMRES applied to \eqref{PBigSchur} as PGMRES in the
sequel.  When $Q_A=A$, it follows from the discussion above that
untruncated AA applied to preconditioned Uzawa is ``essentially
equivalent'' in the above-defined sense to unrestarted PGMRES. When
$Q_B=I$ as well, preconditioned Uzawa reduces to standard Uzawa, and
it follows that untruncated AA applied to standard Uzawa is
``essentially equivalent'' to unrestarted GMRES applied to
\eqref{PBigSchur} with $Q_B=I$, i.e.,
\begin{equation}
 \left(\begin{array}{cc}I & A^{-1}B^T \\[2pt] 0 & \omega
      B{A}^{-1}B^T \end{array}\right)\left( \begin{array}{c}u\\
      p \end{array}\right)=\left( \begin{array}{c}A^{-1}f \\[2pt] \omega
      (B{A}^{-1}f-g)\end{array}\right), \label{BigSchur}
\end{equation}
which mainly requires one solve with $A$ at each iteration.

The forms of \eqref{PBigSchur} and \eqref{BigSchur} are especially
convenient for applying GMRES. For example, the matrix-vector products
required for GMRES applied to \eqref{PBigSchur} mainly require only
one solve with each of $A$ and $Q_B$, in addition to matrix-vector
products with $B$ and $B^T$. Alternatively, one may prefer first to
solve the preconditioned or unpreconditioned Schur-complement system,
i.e.,
\[
Q_{B}^{-1}B{A}^{-1}B^T p = Q_B^{-1} (B{A}^{-1}f-g)\quad {\rm or}\quad
B{A}^{-1}B^T p = B{A}^{-1}f-g, 
\]
and then to form $u=A^{-1}(f-B^Tp)$. This approach has the advantage
of working with shorter vectors in the Arnoldi process within GMRES;
however, it requires an additional solve with $A$ at the end in order
to recover $u$.\footnote{In the important special case in which $A$
  is SPD, one would very likely use the preconditioned
  conjugate-gradient method (PCG) to solve the Schur-complement
  system. Since the approximate solution for $p$ is updated at each
  PCG iteration, one can also update the approximate solution for $u$
  concurrently at very little cost. In this case, the Schur-complement
  approach is likely to be preferred.} In any case, the preconditioned
Schur-complement system points to the need for $Q_B$ to be a good
preconditioner for $BA^{-1}B^T$. See Sections~5 and 10 of
\cite{Benzi05} for more discussion of the Schur-complement approach
and associated preconditioning considerations.

For additional perspectives that reflect recent developments in
splitting preconditioners for incompressible-flow saddle-point
problems, we consider GMRES preconditioned with the Relaxed
Dimensional Factorization (RDF) preconditioner developed in
\cite{Benzi11b}. With $C=0$, we follow \cite{Benzi11b} and write
\eqref{Saddle} equivalently as
\begin{equation}
\left(\begin{array}{cc}
A & B^T  \\
-B & 0 \end{array} \right)
\left(\begin{array}{c}
{ u}  \\
p \end{array} \right)=
\left(\begin{array}{c}
{ f} \\
-g\end{array} \right).
\label{alternatesaddle}
\end{equation}
It is observed in \cite{Benzi11b} that this simple change may, for
some problems (in particular the Stokes equation considered in
Section~4), result in a much more favorable eigenvalue distribution for
Krylov subspace methods. 

The RDF preconditioner for the system (\ref{alternatesaddle}) is an
improvement over the Dimensional Splitting (DS) preconditioner for
solutions of flows with low viscosities \cite{Benzi11a, Benzi11b}. We
note that the DS preconditioner resembles the Hermitian and
skew-Hermitian splitting (HSS) methods, which are popular in the
literature for solving systems like (\ref{alternatesaddle})
\cite{bai2003hermitian, benzi2009generalization}. Recent developments
on HSS-like methods have proposed additional preconditioners for
solving (\ref{alternatesaddle}) when either $C=0$ or $C\neq0$ to improve
convergence compared to standard HSS methods
\cite{chen2015generalized, salkuyeh2015generalized}.

Assuming two-dimensional flow, we partition the velocity vector $u$
into $x$ and $y$ components $u=(u_x,u_y)^T$ and partition $A$ and $B$
conformingly as 
\[
A=\left(\begin{array}{cc}{A}_1&0\\0&{A}_2\end{array}\right), \quad
B = \left(B_1, B_2\right). 
\]
The RDF preconditioner is 
\begin{equation}
M_\beta=\left(\begin{array}{ccc} A_1 &
      -\frac{1}{\beta}B^T_1B_2 & B_1^T \\ 0&A_2&B_2^T\\ -B_1 &
      -B_2 &\beta I
\end{array}\right),
\label{RDFpre}
\end{equation}
where $\beta$ is a problem-dependent relaxation parameter. It is observed
in \cite{Benzi11b} that $M_\beta$ can be factored as 
\[
M_\beta = 
\begin{pmatrix}
I&0&{{1}\over{\beta}}B_1^T\cr0&I&0\cr0&0&I
\end{pmatrix}
\begin{pmatrix}
{\hat A_1}&0&0\cr0&I&0\cr-B_1&0&I
\end{pmatrix}
\begin{pmatrix}
I&0&0\cr0&{\hat A}_2&B_2^T\cr
0&0&\beta I
\end{pmatrix}
\begin{pmatrix}
I&0&0\cr0&I&0\cr0&-{{1}\over{\beta}}B_2&I
\end{pmatrix}, 
\]
where ${\hat A}_i = A_i+{{1}\over{\beta}}B_i^TB_i$ for $i=1$, 2. Thus,
this preconditioner can be economically implemented. It is shown in
\cite{Benzi11b} to be effective on the problems of interest in 
Section~4. In the following, we refer to GMRES preconditioned with
$M_\beta$ simply as RDF.

\section{Numerical experiments}

In this section, we show experimental results comparing the
performance of the methods discussed above on the Stokes and Oseen
problems in several steady incompressible-flow scenarios in two space
dimensions.

\subsection{Fluid equations}
We consider the steady-state incompressible Navier--Stokes equations 
\begin{eqnarray}
-\nu \Delta {u}+({u}\cdot\nabla){u}+\nabla p&={f},\label{NS1}\\
\nabla\cdot{u}&=0,\label{NS2}
\end{eqnarray}
on a bounded, connected domain $\Omega$ in $\mathbb{R}^2$. In
\eqref{NS1}$-$\eqref{NS2}, $\nu$ is the kinematic viscosity,
$u=(u_x,u_y)^T:\Omega\to\mathbb{R}^2$ is the fluid velocity,
$p:\Omega\to\mathbb{R}$ is the pressure field, and
$f:\Omega\to\mathbb{R}^2$ is a force field. The boundary conditions on
$\partial\Omega$ are expressed as $\mathcal{B}{u}=g$, where
$\mathcal{B}$ is a boundary operator. 
%
%\csr{DO WE EVER SAY WHAT $f$ IS IN OUR EXPERIMENTS? IF $f=0$ ALWAYS, WE CAN
%SAY IT HERE. }
%
%\csb{$f$ IS SET TO BE ZERO BUT IT CHANGES BASED ON THE BOUNDARY CONDITIONS (PROBLEM DEPENDENT).}
 
Our interest is in the steady Stokes and Oseen equations, which are
obtained from \eqref{NS1} in two different ways. These are
given, respectively, by
\begin{align}
-\nu \Delta {u}+\nabla p&={f},\label{Stokes}\\
-\nu \Delta {u}+({v}\cdot\nabla){u}+\nabla p&={f},\label{Oseen}
\end{align}
together with the incompressibility condition \eqref{NS2} and boundary
conditions $\mathcal{B}{u}=g$. In the Stokes equations, the nonlinear
term $(u\cdot\nabla)u$ in \eqref{NS1} is dropped with the assumption
that viscous terms dominate and inertial terms can be neglected. In
the Oseen equations, the term $(u\cdot\nabla)u$ in \eqref{NS1} is
linearized and replaced by $({v}\cdot\nabla){u}$, where $v$ is
known. For example, $v$ may be a previous approximate solution
produced in some way. 

Finite-element discretizations of the Oseen and Stokes equations lead
to saddle-point systems of the form \eqref{Saddle}. We focus on the
stable Q2$-$Q1 finite-element discretization with $C=0$, which satisfies
the Lady$\breve{\mbox{z}}$henskaya-Babu$\breve{\mbox{s}}$ka-Brezzi
inf-sup condition \cite{Benzi05}. This is a sufficient condition for
\eqref{Saddle} to have a unique solution and ensures that the
condition number of the Schur complement $BA^{-1}B^T$ is bounded and
independent of the mesh size \cite{Benzi05}.  The matrix $A$ is SPD
for the Stokes equation and is non-symmetric for the Oseen
equation. In our experiments, we used the IFISS MATLAB
package \cite{Elman07} with the Q2$-$Q1 finite-element discretization to
obtain the necessary matrices and vectors for the discretized
problems. In experiments with the Oseen equations, we took $v$ to be
the fifth Picard iterate produced by IFISS applied to the
Navier--Stokes equations \eqref{NS1}$-$\eqref{NS2}. \csb{The initial approximate solution for Picard iteration is to solve the Stokes equation. In our numerical studies, the initial conditions for the pressure and the velocity are zeros.} 

%\csr{WHAT IS THE INITIAL APPROXIMATE SOLUTION FOR IFISS PICARD ITERATION? }
%
%\csb{THE INITIAL PICARD ITERATE IS COMPUTED USING THE STOKES EQUATION. THAT MEANS THE INITIAL CONDITION FOR PICARD ITERATION IS ZERO.}
%
%\csr{WHAT IS THE INITIAL APPROXIMATE SOLUTION FOR NASU, NAPU, ETC? }
%
%\csb{THE INITIAL CONDITION IS ALWAYS ZERO FOR ALL THE METHODS.}

\subsection{Methods and preconditioners}
\label{MethPre}

The methods we compare and our notation for them are as follows:

\smallskip
\begin{tabular}{r l}
NASU & non-accelerated standard Uzawa (\ref{U1})$-$(\ref{U2}) \\ 
ASU & standard Uzawa accelerated with AA \\
NAPU & non-accelerated preconditioned Uzawa (\ref{PU1})$-$(\ref{PU2}) \\
APU & preconditioned Uzawa accelerated with AA \\
PGMRES & GMRES applied to the preconditioned system \eqref{PBigSchur} \\
RDF & GMRES applied to (\ref{alternatesaddle}) preconditioned with
      $M_\beta$ in \eqref{RDFpre}\\
\end{tabular}
\smallskip

With the accelerated algorithms, we indicate the maximum number of
stored residuals in AA (i.e., the truncation parameter $m$) in
parentheses; e.g., ASU(20) denotes ASU with $m=20$ in AA. Similarly,
we show restart parameters in PGMRES and RDF in parentheses, e.g.,
PGMRES(20). In each experiment, truncation and restart parameters were
chosen to be the same so that per-iteration storage and arithmetic
requirements are comparable for the algorithms involved.

We consider the Uzawa algorithm without preconditioning only in
the first test scenario in Section~\ref{channel}. Thereafter, we
consider only the preconditioned algorithm and fix
the preconditioners in \eqref{PU1}--\eqref{PU2} as follows:
\begin{itemize}
\item Stokes equations: $Q_A=A$ and $Q_B$ is the (tridiagonal)
  pressure mass matrix (cf.\cite{Elman94}).
  \item Oseen equations: $Q_A=A$ and 
  \begin{equation}
  \label{LSCprecond}
  Q_B^{-1}=\left(BM_1^{-1}B^T\right)^{-1}BM_1^{-1}AM_1^{-1}B^T\left(BM_1^{-1}B^T\right)^{-1},
    \end{equation}
    where $M_1$ is the diagonal of the velocity mass matrix (cf.\cite{Elman06}). 
\end{itemize}
The preconditioner in \eqref{LSCprecond} is developed by Elman et
al. in \cite{Elman06} as an approximation of the inverse of the Schur
complement. It is a {\em least-squares commutator} preconditioner
specifically referred to as the {\em scaled BFBt method} in
\cite{Elman06}; see also \cite{elman2007least} and
\cite[\S9.2.3]{elman2014finite}. The matrices in (\ref{LSCprecond}) are
easy to obtain since they are available from the saddle-point
system. Applying this preconditioner involves two Poisson-type solves
and one velocity solve.

Note that PGMRES involves the preconditioner $Q_B$, since it appears
in \eqref{PBigSchur}. In experiments with preconditioned Uzawa, $Q_B$
in \eqref{PBigSchur} is as described above. In experiments with standard
(unpreconditioned) Uzawa, $Q_B=I$, and PGMRES is just GMRES applied to
\eqref{BigSchur}. 

The RDF preconditioner $M_\beta$ in \eqref{RDFpre} depends on the
parameter $\beta$. In the experiments reported here, we used values of
$\beta$ that were determined through auxiliary experiments to
approximately minimize numbers of RDF iterations. These values are
shown in each test case and denoted by $\brdf$.

In all numerical experiments, the stopping criterion was 
\begin{equation}
\frac{||r_k||_2}{||b||_2}\leq 10^{-6},\label{stop}
\end{equation}
where $r_k$ is the $k^{th}$ residual and $b$ is the right hand side
vector in \eqref{Saddle}. In the following, we report on the numbers of
iterations required by the algorithms of interest to satisfy
(\ref{stop}) in the test cases described. Since MATLAB was used for
the experiments, these are perhaps more meaningful than timing data. 

\subsection{Channel flow}\label{channel}
The first test problem is channel flow on the square domain
$\Omega=[-1,1]\times[-1,1]$. In this, a parabolic inflow
$u=(1-y^2,0)^T$ is prescribed on the left side of the channel at
$x=-1$. The top and bottom boundary conditions at $y=\pm 1$ are set to
a no-flow Dirichlet condition $u=(0,0)^T$, and a Dirichlet condition
is used for the outflow at $x=1$. With these boundary conditions, we
consider the Stokes equations \eqref{Stokes}$-$\eqref{NS2} with
viscosity $\nu = 1$, for which the standard Uzawa algorithm is
effective but slow. In our experiments, these equations were
discretized on five uniform grids: $16\times16$, $32\times32$,
$64\times64$, $128\times128$, and $256\times256$. These correspond to
659, 2467, 9539, 37507, and 148739 unknowns, respectively.

Although our main interest is in the preconditioned Uzawa algorithm,
for completeness we first consider the standard (unpreconditioned)
Uzawa algorithm \eqref{U1}--\eqref{U2} applied to the Stokes
equations. Specifically, we compare the numbers of iterations
necessary to satisfy \eqref{stop} for ASU, NASU, and PGMRES. (We
remind the reader that PGMRES in this standard-Uzawa case is GMRES
applied to \eqref{BigSchur}.)  Our preliminary experiments showed
that, on this problem, the performance of ASU is not sensitive to the
truncation parameter $m$ in the AA algorithm. In the experiments
reported here, we somewhat arbitrarily chose $m=20$ and took this to
be the restart parameter in PGMRES as well. Also, following
\cite{Elman94}, we took the relaxation parameter in \eqref{U2} to be
$\omega=2/(\lambda_{m}+\lambda_{M})$, where $\lambda_m$ and
$\lambda_M$ are the minimum and maximum eigenvalues of the Schur
complement $BA^{-1}B^T$, which is SPD since $A$ is SPD.

The results of our experiments are shown in
Table~\ref{table:ExactComp}. One sees that ASU(20) consistently
required far fewer iterations than NASU to satisfy \eqref{stop}; thus
acceleration was very effective in these tests. As expected from the
discussion in Section~\ref{PGMRES}, ASU(20) is roughly (but not
exactly) ``essentially equivalent" to PGMRES(20). Additionally, the
convergence of all methods appears to be mesh-independent (i.e., does
not slow appreciably) as the grid is refined. This property is noted
for NASU on the Stokes problem in \cite[\S8.1]{Benzi05}.

\begin{table}[tbh!]
\caption{Iterations required to satisfy \eqref{stop} for the
  channel-flow Stokes problem with $\nu=1$.}
\centering
\begin{tabular}{c|c|c|cc}
\hline
Grid &  ASU (20) & NASU & PGMRES (20)\\  % inserts table %heading
\hline
$16\times16$&		20& 261&	19&\\
$32\times32$&   	26& 268&	29&\\
$64\times64$&     	26& 228& 	29&\\ 
$128\times128$&   	25& 175&	26& \\ 
$256\times256$&    	22& 119&	25& \\ [1ex]
\hline 
\end{tabular}
\label{table:ExactComp}
\end{table}

\begin{figure}[htb!]
\centering
\textbf{(A)}\hspace{2.3in}\textbf{(B)}\\
\includegraphics*[width=0.49\textwidth]{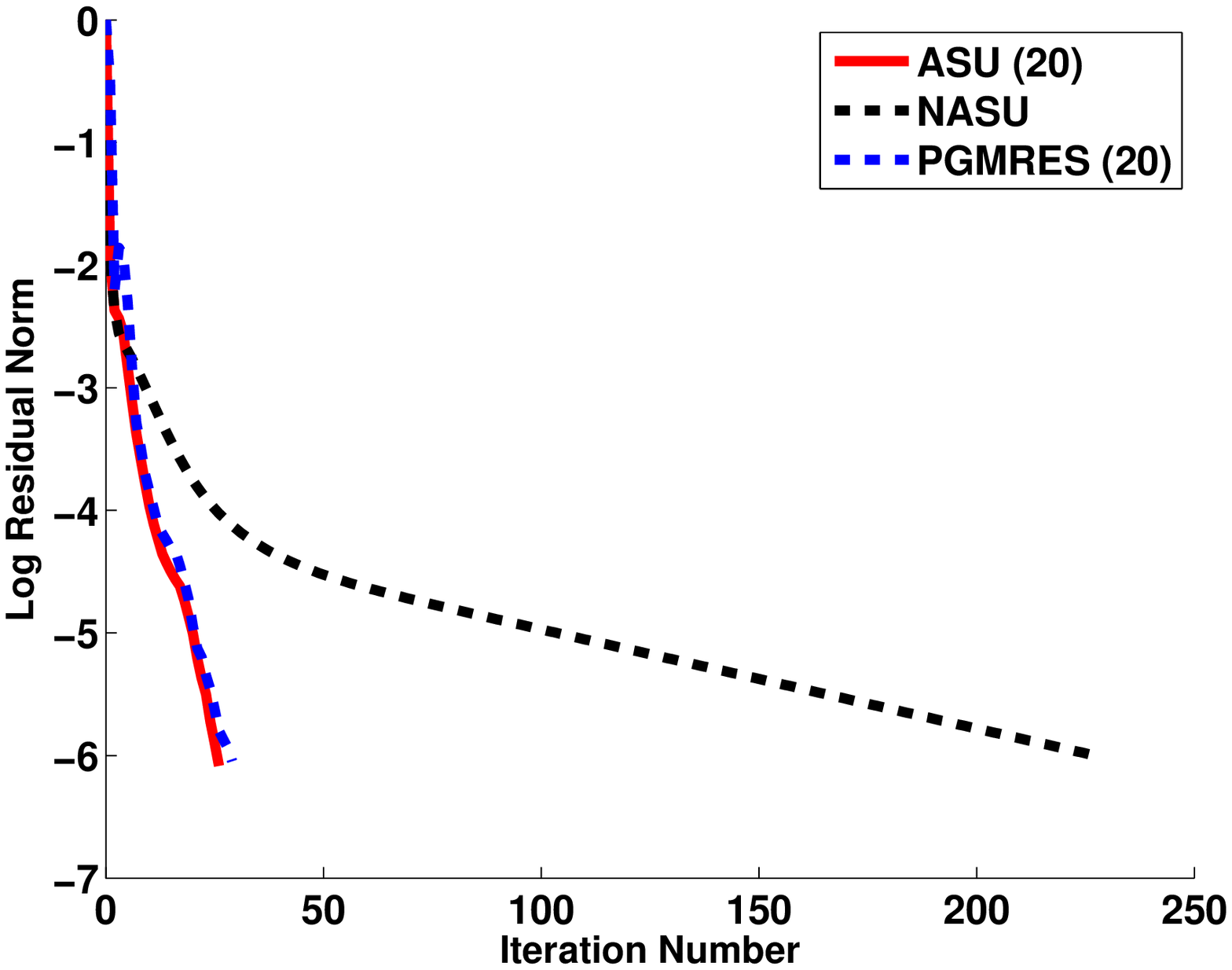}
\includegraphics*[width=0.49\textwidth]{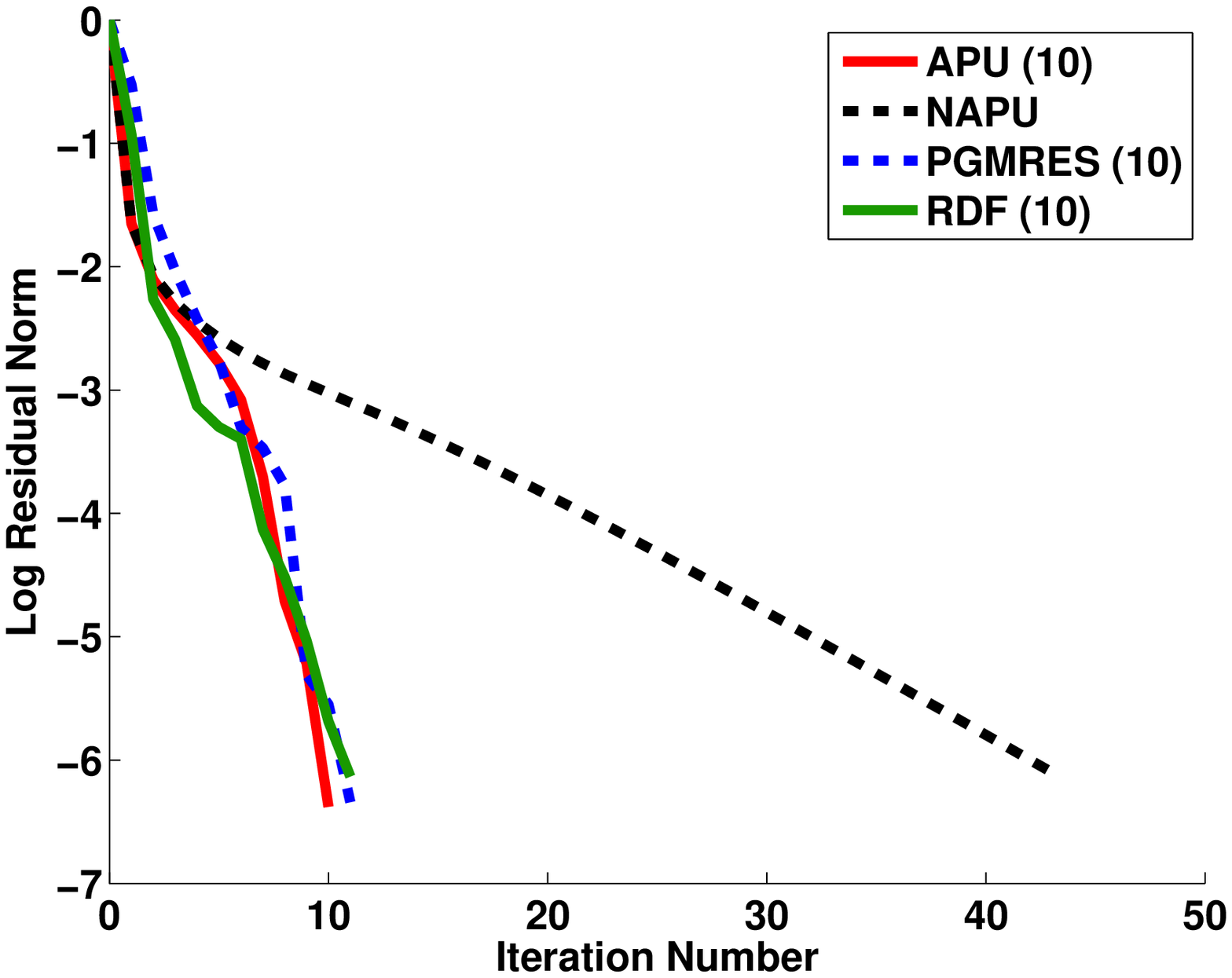}
%\kern-6pt\includegraphics*[width=0.5\textwidth]{FIGURES/AA_PGMRES64x64.eps}
\caption{Log residual-norm plots for the channel-flow Stokes problem
  with: (A) the standard Uzawa algorithm \eqref{U1}--\eqref{U2} for $\nu=1$ on the $64\times64$ grid and (B) the preconditioned Uzawa algorithm \eqref{PU1}--\eqref{PU2} for $\nu=1$ on the $32\times32$ grid. }\label{ResidExact}
\end{figure}

To illustrate the convergence history of the methods, we show in
Fig.~\ref{ResidExact}\,(A) the log residual-norm plots for ASU(20),
NASU, and PGMRES(20) on the 64$\times$64 grid. %In
%Fig.~\ref{ResidExact}\,(B), the log residual-norm plots for ASU(20)
%and PGMRES(20) illustrate the roughly ``essentially equivalent''
%convergence of the two methods.

%\csr{FIG. 1 (B) DOESNT SHOW CONVINCING ESSENTIAL EQUIVALENCE. WHAT TO DO?
%MAYBE OMIT AND JUST SHOW (A)? }
%
%\csb{WE ARE FINE WITH THE IDEA OF REMOVING FIG.~1(B)}

We now consider the preconditioned Uzawa algorithm
\eqref{PU1}--\eqref{PU2} applied to this Stokes problem. The
preconditioners $Q_A$ and $Q_B$ are as described in
Subsection~\ref{MethPre}, and PGMRES is GMRES applied to
\eqref{PBigSchur}. The methods of main interest are APU, NAPU, and
PGMRES. In order to compare with other splitting preconditioners, we
also include results for RDF. As in the unpreconditioned case,
convergence of APU is not sensitive to the truncation parameter $m$ in
AA. Since convergence is somewhat faster than in the unpreconditioned
case, we used $m=10$ in the experiments reported here. We also
observed that convergence is not strongly dependent on the relaxation
parameter $\omega$ in \eqref{PU2} in this case, and we used $\omega=1$
in these experiments.

Table~\ref{channelSt} shows the number of iterations required by the
methods to satisfy \eqref{stop}. As in the unpreconditioned case,
acceleration significantly reduced the necessary number of iterations,
and again the convergence of all methods appears to be mesh-independent.

Figure \ref{ResidExact}\,(B) shows the log residual-norm plots for the
four methods corresponding to the $32\times32$ grid case in
Table~\ref{channelSt}. The log residual-norm plots for APU, PGMRES,
and RDF decrease rapidly. In contrast, the curve for NAPU decreases
much more slowly.
\begin{table}[!htb]
\centering
\caption{Iterations required to satisfy \eqref{stop} for the
  channel-flow Stokes problem with $\nu=1$. In the {\rm RDF} column,
  the numbers in parentheses are the values of $\brdf$ used in each case.}
\begin{tabular}{c|c|c|c|ccccc}
\hline
%\toprule   
%  & \multicolumn{5}{c}{} \\
Grid &APU (10)& NAPU & PGMRES (10)&  RDF(10)    \\
\hline
$16\times16$&    10  & 44 & 10& 10 (0.0044)     \\
$32\times32$   & 10  & 43 & 11& 11 (0.0014)     \\
$64\times64$   & 11  & 41 & 12& 11 (0.0003) \\
$128\times128$&11  & 38 & 12& 10 (0.0001)\\
$256\times256$&11  & 36 & 12& 9 (0.0003)  \\
 \hline
\end{tabular}\label{channelSt}
\end{table}

%\begin{figure}[!htb]
%\centering
%\includegraphics*[width=.5\textwidth]{FIGURES/All_Stokes32x32Channel_New.eps}
%\caption{ Log residual-norm plots for the channel-flow Stokes problem
%  with $\nu=1$ corresponding to the $32\times32$ grid case in
%  Table~\ref{channelSt}. }\label{channelResid}
%\end{figure}

\subsection{Lid-driven cavity}\label{lid}
The second test problem is the ``leaky'' lid-driven cavity problem on
the square domain $\Omega=[-1,1]\times[-1,1]$. With the ``lid'' moving
from left to right, the top boundary condition is $u=(1,0)^T$ at
$y=1$, $0\le x\le1$. No-flow boundary conditions $u=(0,0)^T$ are
imposed on the other three sides. In our experiments, we again
discretized the problem on the five uniform grids used in
Section~\ref{channel}, resulting in the same respective numbers of
unknowns. In these and all subsequent experiments reported here,
results were obtained using only the preconditioned Uzawa algorithm
(\ref{PU1})$-$(\ref{PU2}).

We first consider the Stokes equations \eqref{Stokes}$-$\eqref{NS2} and
take $\nu=1$ in \eqref{Stokes}. In this case, it was sufficient to
use $\omega=1$ in (\ref{PU2}), and $m=10$ led to adequately fast
convergence of APU. The results for this case are summarized in
Table~\ref{StLid}.  As in the channel-flow problem in
Section~\ref{channel}, APU(10), PGMRES(10), and RDF(10) performed
about the same at all grid sizes in our tests, and all three required
significantly fewer iterations to converge than NAPU. Again, the
convergence of all methods appears to be mesh-independent. 

\begin{table}[!htb]
\centering
\caption{Iterations required to satisfy \eqref{stop} for the
  leaky lid-driven cavity Stokes problem with $\nu=1$. In the {\rm
    RDF} column, the numbers in parentheses are the values of $\brdf$
  used in each case.}
\begin{tabular}{c|c|c|c|cccccc}
\hline
%\toprule 
%  & \multicolumn{5}{c}{} \\
Grid& APU(10) & NAPU & PGMRES(10)  & RDF(10)     \\
\hline
$16\times16$   & 12 & 49 & 12 & 12 (0.004)     \\
$32\times32$   & 12  & 50 & 14 & 12 (0.0015)     \\
$64\times64$   & 12  & 50 & 14 & 13 (0.0004) \\
$128\times128$ & 11  & 49 & 14 & 13  (0.0001)   \\
$256\times256$&  11  & 48 & 14 &  13 (0.00002)  &  \\
 \hline
\end{tabular}\label{StLid}
\end{table}

We next consider the Oseen equations \eqref{Oseen}$-$\eqref{NS2} and
report results for $\nu=.1$, $.01$, and $.001$. For this more
challenging problem, we used $m=20$ in APU in all tests. Also, for
$\nu= 0.1$ and $\nu=0.01$, we experimentally determined values of
$\omega$ in \eqref{PU2} that approximately minimized the numbers of
NAPU iterations. For $\nu=.001$, we were unable to find an $\omega$
for which NAPU converged, and so we determined values of $\omega$ that
approximately minimized the numbers of APU iterations (except for the
$16\times16$ grid, for which no method converged for any value of
$\omega$). We did not include RDF in these tests since a complete
account of its performance on this problem is reported in
\cite{Benzi11b}. The results are summarized in Table~\ref{LDCNU}. In
the table, ``$>1000$'' indicates that the iterates appeared to be
converging but failed to satisfy \eqref{stop} within $1000$
iterations; an asterisk indicates that the iterates did not appear to
be converging.

\begin{table}[!htb]
\renewcommand{\arraystretch}{1.1}
\centering
\caption{Iterations required to satisfy \eqref{stop} for the leaky
  lid-driven cavity Oseen problem with $\nu=.1$, $.01$, and $.001$. 
  The values of $\omega$ in \eqref{PU2} shown in the table were
  experimentally determined to approximately minimize the numbers of
  {\rm NAPU} iterations when $\nu=.1$ and $\nu=.01$ and to
  approximately minimize the numbers of {\rm APU} iterations when
  $\nu=.001$.}   
\begin{tabular}{c|c|c|c|c|cccc}
%\hline\hline
%\toprule 
%&&&&&\\[-9.5pt]\cline{2-6} 
%& \multicolumn{5}{c}{}\\[-9.5pt]\cline{2-6} \noalign{\vskip1pt}
%& \phantom{$\big|$} Grid & $\omega$&APU(20) & NAPU & PGMRES(20)    \\
%\cline{2-6} \noalign{\vskip2.3pt} \cline{2-6}
\cline{2-6} 
& Grid & $\omega$&APU(20) & NAPU & PGMRES(20)    \\
%\midrule
\hline%\hline
\multirow{5}{*}{$\nu = .1$} & $16\times16$ 	&0.64   & 10  & 11  & 10    \\
& $32\times32$   		&0.45   & 12  & 17  & 12    \\
& $64\times64$   		&0.29   & 15  & 27  & 15    \\
& $128\times128$ 	&0.16   & 18  & 46  & 18    \\
& $256\times256$ 	&0.087 & 28  & 77   &41     \\
 \hline
\multirow{5}{*}{$\nu = .01$} & $16\times16$ &1.2   & 16  & 51  & 16       \\
& $32\times32$   		&0.74   & 21  & 91  & 20    \\
& $64\times64$   		&0.43   & 23  & 148  & 24    \\
& $128\times128$ 	& 0.24 & 31  & 244  & 40        \\
& $256\times256$ 	&0.12  &  32 &402   &48                \\
 \hline
\multirow{5}{*}{$\nu = .001$} & $16\times16$ &$*$   & $*$  & $*$  & $*$       \\
& $32\times32$   		&1.6   & 99  & $*$  & 378    \\
& $64\times64$   		&0.87   & 111  & $*$  & 600    \\
& $128\times128$ 	&0.31 & 99  & $*$  & $>1000 $       \\
& $256\times256$ 	&0.17   &113   &   $*$    & $>1000 $               \\
\hline
 \end{tabular}\label{LDCNU}
\end{table}

%\csr{WHAT DO THE *'s MEAN (AS OPPOSED TO $>1000$)?}
% 
% \csb{``*" MEANS THAT THE METHOD CAN NOT CONVERGE. ``$>1000$" MEANS THAT THE METHOD IS CONVERGED AT SOME ITERATIONS GREATER THAN 1000. WE CAN USE ``*" FOR THE CASE WHERE THE NUMBER OF ITERATION IS GREATER THAN 1000.}
 
One sees from Table~\ref{LDCNU} that, for $\nu=.1$ and $\nu=.01$,
acceleration not only significantly reduced iteration numbers but also
mitigated mesh dependence, which is at most mild for APU(20) but
pronounced for NAPU. Moreover, for $\nu=.001$, acceleration greatly
improved robustness, with APU(20) succeeding except on the
$16\times16$ grid while NAPU failed on all grids. 
% \csb{We note that
%   the asterisk in Table~\ref{LDCNU} means that the method is failed to
%   converge. In the case of $\nu=.001$, PGMRES(20) converges to the
%   solutions on the grids of $128\times128$ and $256\times256$ for more
%   than 1000 iterations. }

A striking aspect of the results in Table~\ref{LDCNU} is the
difference in performance of APU(20) and PGMRES(20) when the number of
iterations is greater than 20. As discussed in Section~\ref{PGMRES},
the ``essential equivalence'' of APU and PGMRES is theoretically
assured only when the acceleration in APU is not truncated and PGMRES
is not restarted. Evidently, the performance differences seen in
Table~\ref{LDCNU} result from truncating APU and restarting PGMRES,
with the latter having more significant negative effects on
performance than the former.

 To further investigate the effects of the choice of the restart value
on PGMRES convergence on this problem, we tried a number of different
restart values with $\nu=.001$ on the $64\times64$ grid. The results
are in Table~\ref{PGMRESrestartLid}. One sees that performance
comparable to that of APU(20) in Table~\ref{LDCNU} was achieved, but
only with a restart value much greater than 20 and a correspondingly
greater cost of storage and arithmetic.

\begin{table}[!htb]
\centering
\caption{Iterations required by restarted {\rm PGMRES} to satisfy \eqref{stop} for the leaky
  lid-driven cavity Oseen problem with $\nu=.001$ on the $64\times64$
  grid with $\omega=0.87$ (cf. Table~\ref{LDCNU}). }  
\begin{tabular}{c|c|c|c|c|c|c}
\hline
Restart Value & 20 & 30 & 40 & 50 & 60 & 70 \\% & 80 \\
\hline
Iterations & 600 & 343 & 190 & 167 & 131 & 98 \\ % & 102 \\
 \hline
\end{tabular}  \label{PGMRESrestartLid}
\end{table}

\begin{figure}[!htb]
\centering
\textbf{(A)}\hspace{2in}\textbf{(B)}\\
\includegraphics*[width=.5\textwidth]{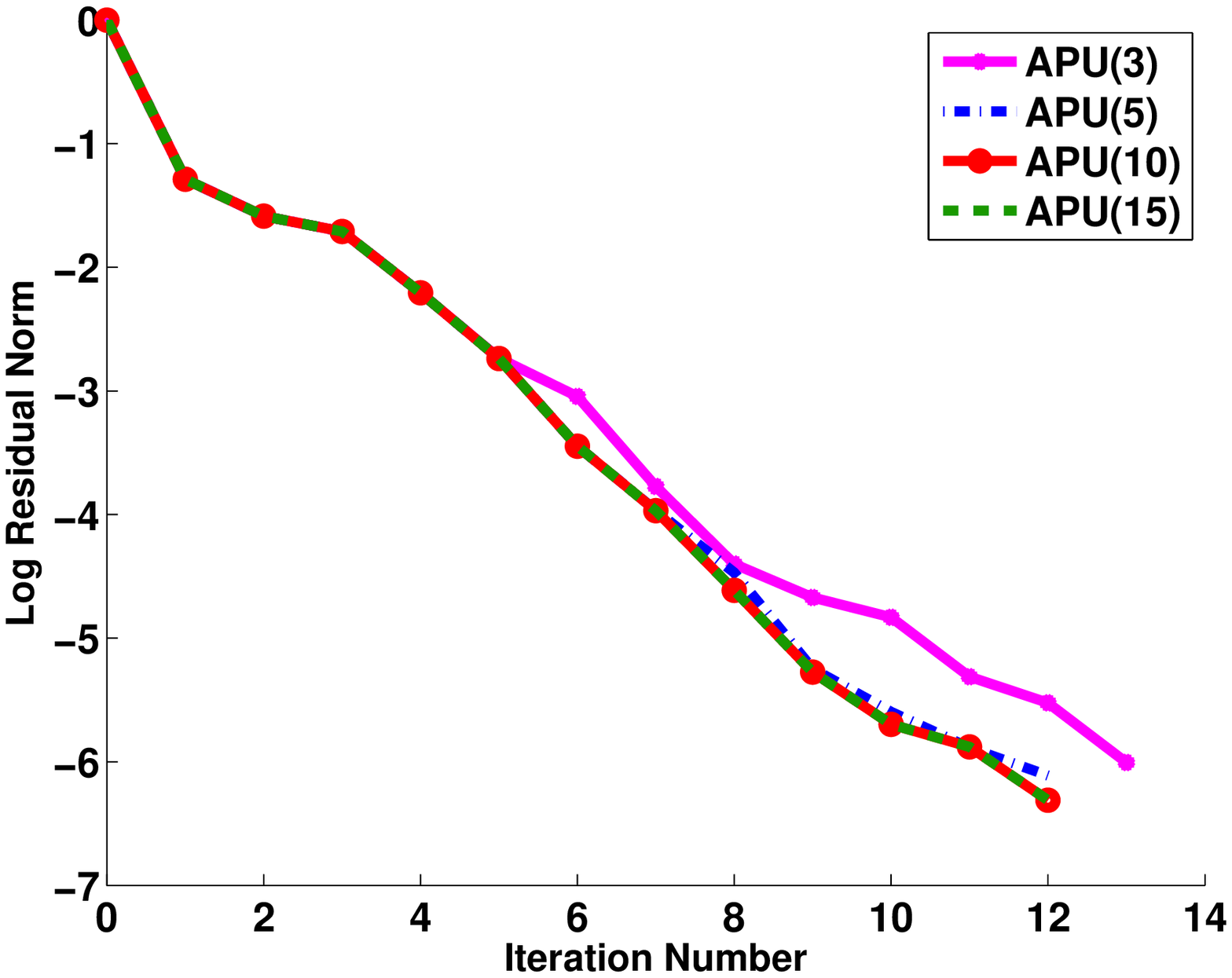}
\kern-6pt\includegraphics*[width=.5\textwidth]{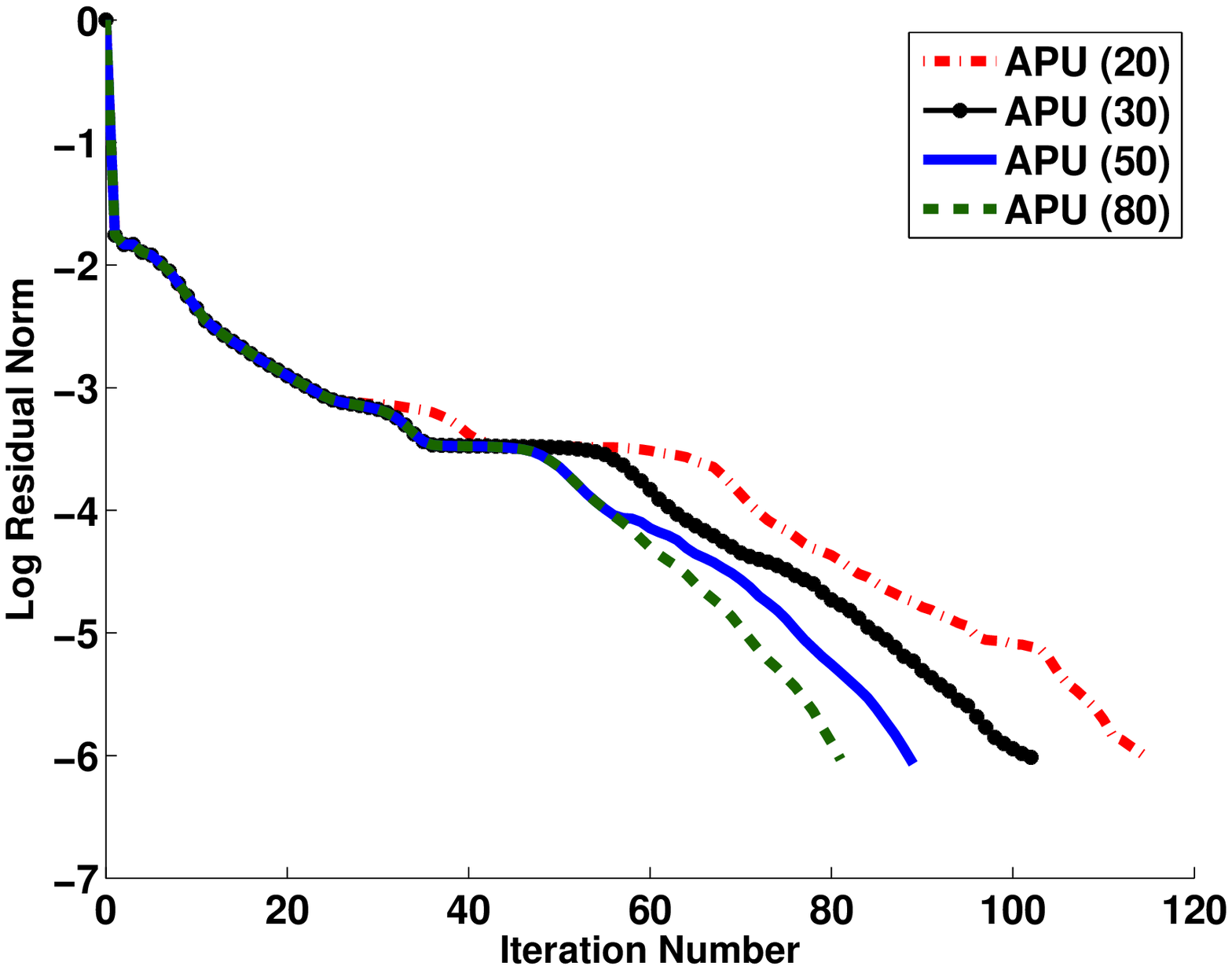}
\caption{Log residual norm plots for the leaky lid-driven cavity
  problem using different numbers of stored residuals $m$ in {\rm APU}
  on the $64\times64$ grid: {\rm (A)} for the Stokes equations with
  $\nu=1$ in {\rm \eqref{Stokes}} and $\omega=1$ in {\rm \eqref{PU2}}
  (cf. Table~~\ref{StLid}); {\rm (B)} the Oseen equations with
  $\nu=0.001$ in {\rm \eqref{Stokes}} and $\omega=0.87$ in {\rm
    \eqref{PU2}} (cf. Table~\ref{LDCNU}).}
\label{varym}
\end{figure}

We also explored the sensitivity of APU convergence to the parameter
$m$, the number of stored residuals in AA. In Fig.~\ref{varym}, we
illustrate how the convergence of APU varies with $m$ for the Stokes
and Oseen equations. For the Stokes equations, the log residual-norm
plots decrease smoothly and show relatively little sensitivity to
$m$, with little meaningful difference in the plots for the four
values of $m$. For the Oseen equations, the plots decrease smoothly
for the most part, although there is a modest plateau region in the
middle iterations. There is also somewhat greater sensitivity to $m$
in this more challenging case.

%\csr{THE RESULTS IN FIG. 3(B) FOR APU(20) DON'T AGREE WITH
%TABLE~\ref{LDCNU}, WHICH SHOWS 111 ITERATIONS REQUIRED BY APU(20).}
%
%\csb{FIG. 3(B) IS UPDATED TO REFLECT THE CORRECT NUMBER OF ITERATIONS FOR APU(20). FIG. 3(A) IS ALSO GOT UPDATED WITH BIGGER FONT.}

\subsection{Flow past an obstacle}\label{obstacle}
The next test problem is flow past an obstacle. This corresponds to
channel flow inside the domain $\Omega=[0, 8]\times[-1,1]$ with a
square object located at $[0,2]\times[0.5,0.5]$. The boundary
conditions are as in Subsection~\ref{channel} \csb{except Neumann boundary conditions are set for the outflow}. We consider four
different grid sizes $16\times32$, $32\times64$, $64\times128$, and
$128\times256$ corresponding to 2488, 9512, 37168, and 146912
unknowns, respectively. We again focus on the preconditioned Uzawa
algorithm \eqref{PU1}$-$\eqref{PU2} with and without acceleration
applied to the Oseen equations \eqref{Oseen}$-$\eqref{NS2} and take
$m=20$ in APU($m$). We include results for RDF, since to our knowledge
these have not been reported elsewhere. For this problem and the
backward-facing step problem considered next, the smallest value of
$\nu$ considered was $.005$, rather than $.001$, consistent with
observations in \cite{Benzi11b} that a steady solution would be
unstable for the latter problem with $\nu=.001$. The results are
summarized in Table~\ref{FPONU}.

%\csr{ IS WHAT I SAID ABOUT THE BOUNDARY CONDITIONS CORRECT? }

\begin{table}[!htb]
\centering
\caption{Iterations required to satisfy \eqref{stop} for the
  obstacle-flow Oseen problem with $\nu=.1$, $.01$, and $.005$. The
  values of $\omega$ in \eqref{PU2} shown in the table were
  experimentally determined to approximately minimize the numbers of
  {\rm NAPU} iterations for each value of $\nu$. In the {\rm RDF}
  column, the numbers in parentheses are the values of $\brdf$ used in
  each case.}  

\begin{tabular}{c|c|c|c|c|c|c}
\cline{2-7} 
& Grid& $\omega$&APU(20) & NAPU & PGMRES(20) & RDF(20)   \\
%\midrule
\hline
\multirow{4}{*}{$\nu = .1$} & $16\times32$ 	&0.72   & 12  & 13 & 11  &  16 (0.044)   \\
& $32\times64$   	&0.49   & 13  & 17 & 13  &  17 (0.014)\\
& $64\times128$ 	&0.34   & 16  & 25 & 15  &  17 (0.005)\\
& $128\times256$ &0.21  & 19  & 39 & 18 & 16 (0.001)     \\
 \hline
\multirow{4}{*}{$\nu = .01$} & $16\times32$&1.61&34& 138&36&18 (0.096)   \\
& $32\times64$   	&1.12   & 25  & 67  & 26  &  20 (0.041)\\
& $64\times128$ 	&0.85   & 16  & 21  & 16  &  22 (0.013)\\
& $128\times256$ &0.49 &  14 & {19}   &  16 &  24 (0.004)    \\
 \hline
\multirow{4}{*}{$\nu = .005$} & $16\times32$&1.28&55&489& {66}&21 (0.100)   \\
& $32\times64$   	&0.76   & {50}  & 332  & 57  &  {19} (0.050)\\
& $64\times128$ 	&0.76   & 28  & 77  & 30  &  31 (0.013)\\
& $128\times256$ &0.52 &  16  & 36&  16&   34 (0.005)   \\
 \hline
\end{tabular}\label{FPONU}
\end{table}

The results in Table~\ref{FPONU} show that, as in the previous
test problems, acceleration was effective in reducing the number of
preconditioned Uzawa iterations. When $\nu=.1$, acceleration also
mitigated mesh dependence somewhat, although it is fairly mild for
NAPU in this case. Interestingly, when $\nu=.01$ and $\nu=.005$, the
iteration numbers required by APU(20), NAPU, and PGMRES(20) all
decrease as the grid is refined. For all $\nu$ values, the iteration
numbers for APU(20) and PGMRES(20) consistently differ by little, in
contrast to the previous problem. The iteration numbers for RDF(20)
suggest at most mild mesh dependence for that method. They are about
the same as those of APU(20) and PGMRES(20) when $\nu=.1$. For the two
smaller values of $\nu$, they are somewhat smaller for the two coarser
grids and somewhat larger for the two finer grids.

\begin{figure}[!htb]
\centering
\textbf{(A)}\hspace{2in}\textbf{(B)}\\
%\textbf{(A)}
%\hskip-34pt
\includegraphics*[width=.5\textwidth]{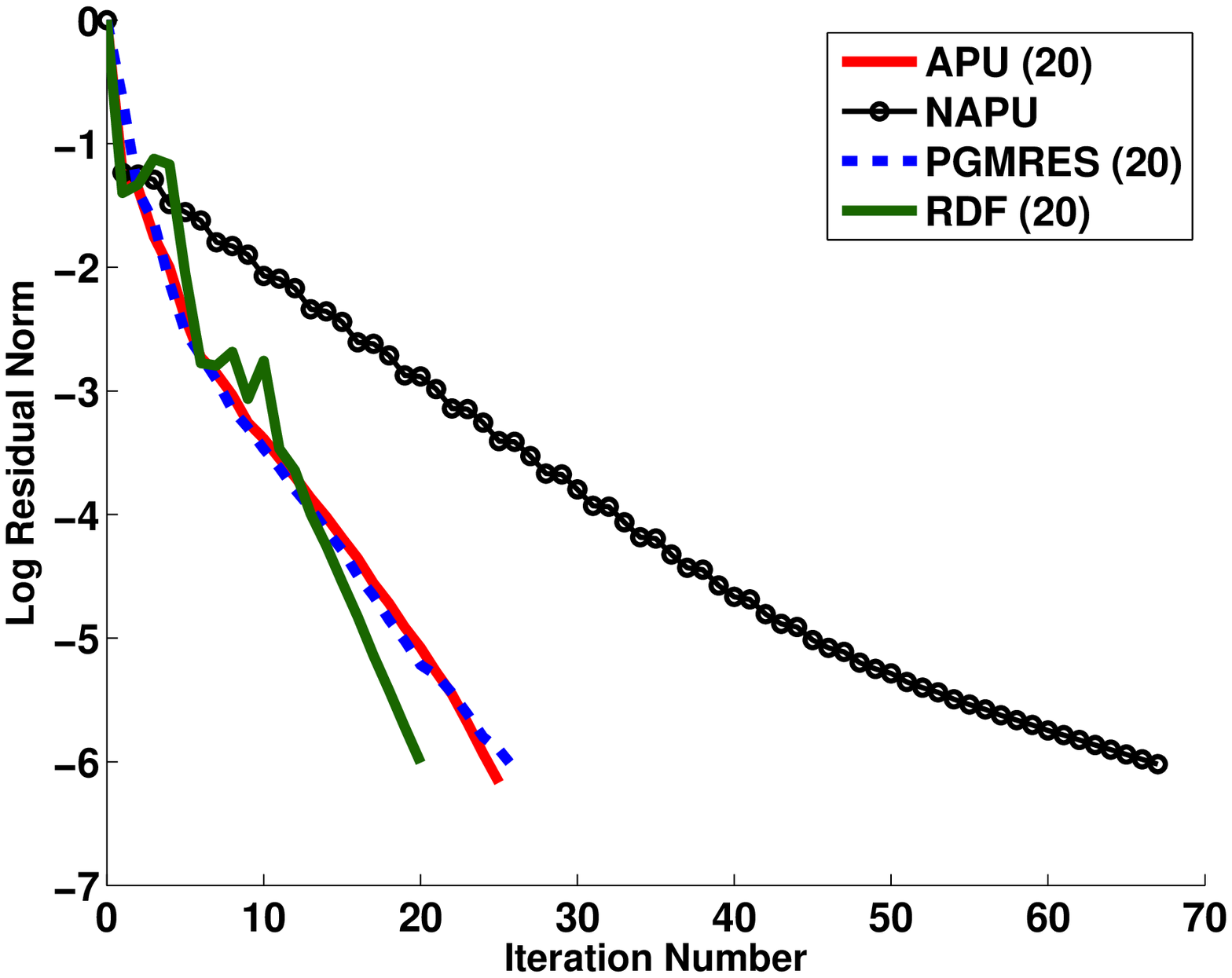}
\hskip-6pt\includegraphics*[width=.5\textwidth]{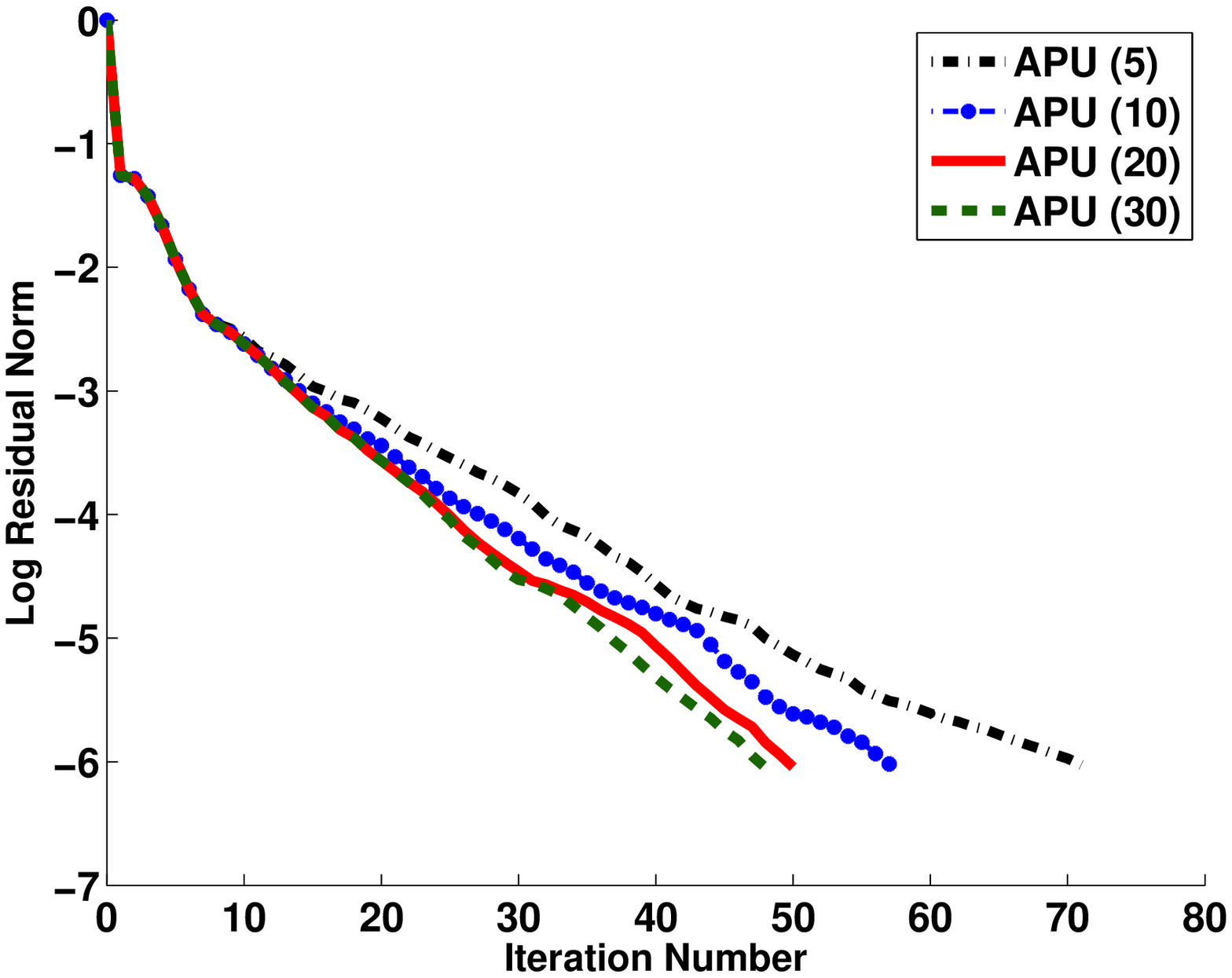}
\caption{{Log residual-norm plots for the obstacle-flow Oseen problem
    on the $32\times64$ grid: (A)~all methods with $\nu=.01$; (B)~{\rm
      APU($m$)} for different $m$ with $\nu=.005$
    (cf. Table~\ref{FPONU}). }}  
\label{flowpastobs}
\end{figure}

To further illustrate the convergence of the methods, we show in
Fig.~\ref{flowpastobs}\,(A) the log residual-norm plots of APU(20),
NAPU, PGMRES(20), and RDF(20) for the case $\nu=0.01$ and the
$32\times64$ grid. Also, we note that the choice of $m$, the number of
stored residuals in AA, does matter for this Oseen problem. In
Fig.~\ref{flowpastobs}\,(B), we show results for different choices of
$m$ for this problem with $\nu=0.005$ on the $32\times64$ grid. The
figure clearly indicates that APU($m$) converges faster as $m$
increases. It is also notable that significant improvement over NAPU
was obtained with $m$ as small as five. 

%\csr{THE RESULTS IN FIG. 4(B) FOR APU(20) DON'T AGREE WITH
%TABLE~\ref{FPONU}, WHICH SHOWS 55 ITERATIONS REQUIRED BY
%APU(20). }
%
%\csb{THE OLD FIGURE IS SWAPPED OUT WITH THE NEW ONE WITH CORRECT CAPTION.}
%
%\csr{ALSO, IN FIG. 4(A), THE CAPTION SHOULD SAY APU(20), PGMRES(20) AND
%RDF(20).}
%
%\csb{THE FIGURE IS CORRECT (I CHANGED THE FONT TO MAKE IT LOOK BETTER). THE TABLE IS UPDATED TO REFLECT THE CORRECT NUMBER OF ITERATIONS.}
 
\subsection{Backward-facing step}\label{step}
Our final test problem is flow over a backward-facing step. In this,
the domain is $\Omega=[-1,0]\times[0,1]\cup[0,5]\times[-1,1]$, which
has a downward ``step'' at $x=0$.  As in Subsection~\ref{channel}, a
parabolic flow profile is imposed at the inflow boundary ($x=-1$ and
$0\leq y\leq 1$). Neumann boundary conditions are imposed at the
outflow boundary ($x=5$ and $-1< y< 1$).  At the other walls, a
no-flow boundary is prescribed. We study the problem using four
different grids having the general form $2^n\times3\cdot2^n$ for
$n=4, 5, 6, 7$, corresponding to 1747, 6659, 25987, and 102659
unknowns, respectively. As in Subsection~\ref{obstacle}, we consider
the Oseen equations \eqref{Oseen}$-$\eqref{NS2} with $\nu=.1$,
$\nu=.01$, and $\nu=.005$ and report on the performance of the methods
with $m=20$ in APU($m$). The results are summarized in
Table~\ref{BFSNU}.

\begin{table}[!htb]
\centering
\caption{Iterations required to satisfy \eqref{stop} for the
  backward-facing step Oseen problem with $\nu=.1$, $.01$, and
  $.005$. The values of $\omega$ in \eqref{PU2} shown in the table were
  experimentally determined to approximately minimize the numbers of
  {\rm NAPU} iterations for each value of $\nu$. In the {\rm RDF}
  column, the numbers in parentheses are the values of $\brdf$ used in
  each case.} 

\begin{tabular}{c|c|c|c|c|c|c}
\cline{2-7} 
& Grid& $\omega$&APU(20) & NAPU & PGMRES(20) & RDF(20)   \\
%\midrule
\hline
\multirow{4}{*}{$\nu = .1$} & $16\times48$ 	&0.64   & 11  & 12  & 11  &  15 (0.035)   \\
& $32\times96$   	&0.45   & 14  & 17  & 14  &  14 (0.014)\\
& $64\times192$ 	&0.29   & 17  & 26  & 17  &  13 (0.005)\\
& $128\times384$ &0.17 & 20& 40   &  19 &  14 (0.0009)    \\
\hline
\multirow{4}{*}{$\nu = .01$} & $16\times48$ &0.93   & 24  & 61  & 24  &  15 (0.140)   \\
& $32\times96$   	&0.55   &{15}  & 25  & 14  &  15 (0.043)\\
& $64\times192$ 	&0.46   & 14  & 24  & 14  &  16 (0.013)\\
& $128\times384$ &0.38 &   19& 32   & 20	&  17 (0.005)    \\
\hline
\multirow{4}{*}{$\nu = .005$} & $16\times48$ &0.44   & 41  & 400  & 40  &  20 (0.170)   \\
& $32\times96$   	&0.42   & 28  & 112  & 31  &  {19 (0.061)}\\
& $64\times192$ 	&0.32   & 20  & 49  & 20  &  19 (0.017)\\
& $128\times384$ &0.31&   21& 36   &   20	& 19 (0.006)     \\
\hline
\end{tabular}\label{BFSNU}
\end{table}

\begin{figure}[!htb]
\centering
\includegraphics*[width=.45\textwidth]{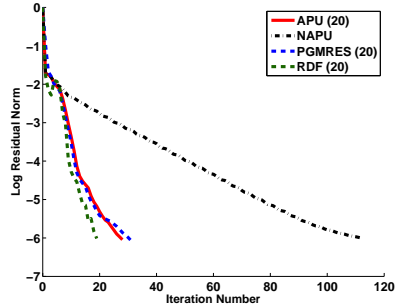}
\caption{{Log residual-norm plots for all methods on the
    backward-facing step Oseen problem with viscosity $\nu=0.005$ and
    the $32\times96$ grid
    (cf. Table~\ref{BFSNU}). }}\label{resBFS0005}  
\end{figure}

The results in Table~\ref{BFSNU} are qualitatively very similar to
those in Table~\ref{FPONU}. As in that case, acceleration was
effective in reducing the number of preconditioned Uzawa iterations
and also in mitigating the rather mild mesh dependence of NAPU when
$\nu=.1$. When $\nu = .01$ and $\nu = .005$, the table shows some
tendency for the iteration numbers of APU(20), NAPU, and PGMRES(20) to
actually decrease as the grid is refined, although it is not as
pronounced or consistent as in the previous case. The iteration
numbers for APU(20) and PGMRES(20) are very similar and do not differ
much from those of RDF(20), which show no significant mesh
dependence. We note that this problem is considered in \cite{Benzi05}
but only for $\nu=.005$ on the $16\times48$, $32\times96$, and
$64\times192$ grids. In these cases, the results for RDF(20) shown
here are slightly more favorable than those in \cite{Benzi05},
presumably because of slightly different values of $\brdf$ (denoted by
$\alpha$ in \cite{Benzi05}). We conclude this section by showing in
Figure~\ref{resBFS0005} the performance of all methods in the case of
$\nu=.005$ and the $32\times96$ grid.

%
%\csr{IN FIG. 5, THE CAPTION SHOULD SAY APU(20), PGMRES(20) AND RDF(20).}
%
%\csb{THE OLD FIGURE IS SWAPPED OUT WITH THE NEW ONE WITH CORRECT CAPTION.}

\section{Discussion \& Conclusions}
We have introduced {\em Anderson acceleration} (AA) as an easily
implemented and economical way of improving the convergence of the
Uzawa algorithm for solving a saddle-point system
\eqref{Saddle}. While such systems arise in many settings, the
applications of primary interest here are the steady Stokes and Oseen
problems of incompressible flow. As far as we know, this is
the first time that AA has been used either with the Uzawa algorithm
in any context or with any iterative method for these Stokes and Oseen
problems.

After introducing the standard (unpreconditioned) and preconditioned
Uzawa algorithms in Section~\ref{varUzawa}, we describe AA in
Section~\ref{accUzawa} and discuss its application to the Uzawa
algorithm. We show that, viewed as a fixed-point iteration, the
preconditioned Uzawa algorithm is a stationary iteration for
\eqref{Saddle} determined by the splitting \eqref{Split}. It follows
from results in \cite{Walker11} that preconditioned Uzawa accelerated
with untruncated AA is ``essentially equivalent'' in a certain sense
to unrestarted GMRES applied to the preconditioned system
\eqref{LPSystem}. When the Uzawa preconditioning matrix $Q_A$ in
\eqref{PU1} is equal to $A$, the preconditioned system
\eqref{LPSystem} simplifies to \eqref{PBigSchur}, a form especially
convenient for applying GMRES. Similar results hold for standard Uzawa
as a special case, with the simplified form of the preconditioned
system given by \eqref{BigSchur}.

In Section~4, we present results of extensive numerical experiments
comparing the performance of the Uzawa algorithm with and without
acceleration in several steady Stokes and Oseen incompressible-flow
scenarios. Our main interest is in the preconditioned Uzawa algorithm
\eqref{PU1}$-$ß\eqref{PU2}, with preconditioning as described in
Section~\ref{MethPre}; however, we also report results for the
unpreconditioned algorithm \eqref{U1}$-$\eqref{U2} in one case, the
channel-flow Stokes problem considered in Section~\ref{channel}. In
discussing the test results below, we often refer to methods using the
acronyms given in Section~\ref{MethPre}; for example, PGMRES refers to
GMRES applied to the preconditioned system \eqref{PBigSchur}.

Although, as noted above, unrestarted PGMRES is ``essentially
equivalent'' to Uzawa accelerated with untruncated AA, this
equivalence does not hold when PGMRES is restarted and AA is
truncated. Consequently, we included restarted PGMRES in our tests in
order to see how its performance compared to that of Uzawa accelerated
with truncated AA. To broaden the comparison, we also included in most
test cases GMRES applied to the equivalent system \eqref{RDFpre}
preconditioned with the sophisticated Relaxed Dimensional
Factorization preconditioner developed in \cite{Benzi11b}, referred to
simply as RDF here.

In all test cases, acceleration reduced the numbers of iterations
required by the Uzawa algorithm to satisfy the convergence criterion
\eqref{stop}. The reduction was usually significant, and in many cases
it was dramatic. Accelerated Uzawa usually required about the same
numbers of iterations or modestly fewer than restarted PGMRES;
however, in one case, that of the leaky lid-driven cavity problem (see
Table~\ref{LDCNU}), the accelerated method APU(20) was significantly
more robust than PGMRES(20) and often required significantly fewer
iterations when the iterates from both methods converged. The results
in Table~\ref{PGMRESrestartLid} suggest that the PGMRES iterates can
be made to converge in as few iterations as those of APU(20) in this
test case, but only by taking the restart value to be much larger than
20. 

Accelerated Uzawa usually required about the same numbers of
iterations as RDF, although not in every case (see, e.g.,
Table~\ref{FPONU}). We note that the RDF preconditioner $M_\beta$ in
\eqref{RDFpre} depends on a parameter $\beta$. Although some
guidelines for choosing $\beta$ are given in \cite{Benzi11b},
determining effective values of $\beta$ may require some effort in
practice. In our tests, values of $\beta$ were determined through
auxiliary experimentation to give near-optimal convergence in each
case. 

We also used auxiliary experimentation to determine values of the
parameter $\omega$ in the standard and preconditioned Uzawa iterations
\eqref{U1}$-$\eqref{U2} and \eqref{PU1}$-$\eqref{PU2} that approximately
minimized the numbers of non-accelerated Uzawa iterations in each test
case. Developing more efficient methods for determining
effective choices of $\beta$ in \eqref{RDFpre} and $\omega$ in the
Uzawa iteration will be a subject of future work. 

In the tests involving Stokes flow (see
Tables~\ref{table:ExactComp}-\ref{StLid}), all methods showed little
adverse mesh dependence. This is not surprising, in view of the
mesh-independence property of non-accelerated Uzawa for Stokes flow
noted in \cite[\S8.1]{Benzi05}. In the tests with Oseen flow (see
Tables~\ref{LDCNU}-\ref{BFSNU})), the RDF method showed at most very
modest mesh dependence in all cases. With the other methods, mesh
dependence varied considerably. For the leaky lid-driven cavity Oseen
problem (Table~\ref{LDCNU}), mesh dependence of the non-accelerated
preconditioned Uzawa algorithm NAPU was pronounced, while that of the
accelerated method APU(20) was much more benign. For the obstacle-flow
and backward-facing step Oseen problems (Tables~\ref{FPONU} and
\ref{BFSNU}), convergence of both NAPU and APU(20) was not very
adversely affected and, in some cases, actually improved as the mesh
was refined.

% Not sure about including following two observations. 

% Specifically, APU appeared to be mesh-independent and was a
% competitive scheme compared to Krylov subspaces methods for the Stokes
% equations.  For the steady Oseen problems, APU was not
% mesh-independent for lower viscosity and larger grid sizes. This is
% partially due to the Schur complement appearing in our fixed-point
% method, as seen in \ref{StokesU} \cite{Benzi05}.

% If $\omega$ is not chosen carefully, this could lead to
% ill-conditioning of the system. 

% Following moved to Section 3.3: 

% The preconditioner of the RDF method for
% the system (\ref{alternatesaddle}) is an improvement over the
% Dimensional Splitting (DS) preconditioner for solutions of flows with
% low viscosities \cite{Benzi11a, Benzi11b}. We note that the DS
% preconditioner resembles the Hermitian and skew-Hermitian splitting
% (HSS) methods. This class of preconditioners are popular in the
% literature for solving systems like (\ref{alternatesaddle})
% \cite{bai2003hermitian, benzi2009generalization}. Recent developments
% on HSS-like methods have proposed different preconditioners for
% solving (\ref{alternatesaddle}) when $C=0$ and $C\neq0$ to accelerate
% the convergence of the solution when compared to the standard HSS
% methods \cite{chen2015generalized, salkuyeh2015generalized}.

There are other Uzawa-like methods for iteratively solving the
saddle-point system (\ref{Saddle}) that we have not considered
here. For instance, if $Q_A^{-1}$ in (2.2a) is replaced by
$Q_{A}^{-1}=\eta P^{-1}$, where $\eta$ is a positive relaxation
parameter and $P$ is a preconditioner, then the system
(\ref{PU1})--(\ref{PU2}) is called the parameterized inexact Uzawa
(PIU) method \cite{bai2005generalized}. Several frameworks have been
developed for the PIU method to accelerate convergence to solutions
for both the symmetric and nonsymmetric generalized saddle-point
problems \cite{bai2008parameterized}. Another class of PIU methods has
been developed recently to iteratively solve for solutions of complex
symmetric linear systems \cite{zheng2015fast}; these methods accelerate
convergence using a correction technique. The corrected PIU methods are
also shown to have faster convergence compared to some Uzawa-type and
Hermitian or skew-Hermitian splitting-type methods
\cite{zheng2015fast}. Other Uzawa-like methods have been proposed in
recent years to solve saddle-point systems. These methods include
Uzawa-SOR, Uzawa-AOR, Uzawa-SAOR, and Uzawa-SSOR \cite{zhang2010class,
  yun2013variants}. These methods utilize the decomposition of the
matrix $A$ as $A=D-L-U$, where $D, L,$ and $U$ are diagonal, strictly
lower-triangular, and strictly upper-triangular, respectively. The
proposed methods are effective in solving the saddle-point system and
have less computational time compared to the generalized SOR method
\cite{zhang2010class}.

While we have focused on fluid flow governed by the Stokes and Oseen
equations, Uzawa-like methods have also been applied to other fluid
models, including non-Newtonian fluids such as viscoplastic and
viscoelastic fluids. For instance, a convergence study and numerical
simulations for the Bingham model of a viscoplastic fluid has been
investigated in \cite{dean2007numerical}. Simulations for unsteady
Bingham flow in cylinders and in a lid-driven cavity using the
Newton-Conjugate Gradient-Uzawa algorithm were studied where the problems were
discretized using finite-element methods \cite{dean2007numerical}.
% A recent study on Bingham flow has done simulations using a
% finite-difference method and the computational experiments were
% carried out for high Reynolds number\cite{muravleva2015uzawa}.
For viscoelastic fluids, a time-dependent Oldroyd model has been
treated using an Uzawa-like scheme called the gauge-Uzawa
(projection-Uzawa) finite-element method \cite{si2015gauge}. In this
method, the velocity component is decomposed into an unknown vector
and an auxiliary scalar variable, easily handling boundary
derivatives. This method, however, does not require the formation of
the saddle-point system. In the future, we will explore extending the
accelerated Uzawa schemes to solve non-Newtonian fluid problems.

\bibliographystyle{siam}
\bibliography{BibCopperMountain}
\end{document}